\documentclass[12pt]{article}
\usepackage{amsmath,amsfonts,latexsym,amsthm,amssymb}
\newcommand{\RR}{\mathbb R}
\newcommand{\Qp}{\mathbb Q_p}
\newcommand{\D}{D^\alpha}
\newcommand{\II}{I^\alpha}

\numberwithin{equation}{section}

\DeclareMathOperator{\uu}{\overline{u}}
\begin{document}

\newtheorem{lem}{Lemma}
\newtheorem{teo}{Theorem}
\newtheorem*{cor}{Corollary}
\newtheorem*{prop}{Proposition}

\pagestyle{plain}
\title{Nonlinear Pseudo-Differential Equations for Radial Real Functions on a Non-Archimedean Field}
\author{Anatoly N. Kochubei
\\ \footnotesize Institute of Mathematics,\\
\footnotesize National Academy of Sciences of Ukraine,\\
\footnotesize Tereshchenkivska 3, Kiev, 01024 Ukraine\\
\footnotesize E-mail: \ kochubei@imath.kiev.ua}
\date{}
\maketitle

\vspace*{3cm}
\begin{abstract}
In an earlier paper (A. N. Kochubei, {\it Pacif. J. Math.} 269 (2014), 355--369), the author considered a restriction of Vladimirov's fractional differentiation operator $D^\alpha$, $\alpha >0$, to radial functions on a non-Archimedean field. In particular, it was found to possess such a right inverse $I^\alpha$ that the change of an unknown function $u=I^\alpha v$ reduces the Cauchy problem for a linear equation with $D^\alpha$ (for radial functions) to an integral equation whose properties resemble those of classical Volterra equations. In other words, we found, in the framework of non-Archimedean pseudo-differential operators, a counterpart of ordinary differential equations. In the present paper, we study nonlinear equations of this kind, find conditions of their local and global solvability.
\end{abstract}

\vspace{2cm}


\newpage
\section{Introduction}

The basic linear operator defined on real- or complex-valued functions on a non-Archimedean local field, for example the field $\Qp$ of $p$-adic numbers, is the Vladimirov pseudo-differential operator $\D$ of fractional differentiation; see \cite{K2001,VVZ}; for further development of this and related subjects see \cite{AKS,KKZ,ZG}.

It was found in \cite{K2014} that properties of $\D$ become much simpler on radial functions. For this case, it resembles the Caputo-Djrbashian fractional differentiation of real analysis. In particular, this operator possesses a right inverse, which can be seen as a $p$-adic counterpart of the Riemann-Liouville fractional integral or, for $\alpha =1$, the classical anti-derivative. As an application, in \cite{K2014} we considered linear equations for radial functions, in which the restriction of $\D$ plays the role of a derivative.

In this paper, we extend the above framework to nonlinear equations. We consider equations of the form
\begin{equation}
\label{1.1}
\left( \D u\right) (|t|_K) =f(|t|_K,u(|t|_K)),\quad 0\ne t\in K,
\end{equation}
with the initial condition
\begin{equation}
\label{1.2}
u(0)=u_0
\end{equation}
where $K$ is non-Archimedean local field with the normalized absolute value $|\cdot |_K$ (see section 2 for the definitions).

Under certain Lipschitz type condition, we prove the existence and uniqueness of a local solution of the Cauchy problem (\ref{1.1})-(\ref{1.2}) defined on a sufficiently small ``interval'' of $t$. While this property is similar to classical ones, our analytic technique is different, due to the different nature of our $p$-adic fractional integral. The analysis of the structure of (\ref{1.1}), that is of the properties of $\D$ on radial functions, then shows that under some assumptions regarding the function $f$, the local solution is extended to a global one defined on the whole field $K$. While for linear equations \cite{K2014} we used the theory of compact operators, for nonlinear equations studied in this paper, the above results are obtained by direct constructions based on iteration processes.

\section{Preliminaries}

{\bf 2.1. Local fields.} Let $K$ be a non-Archimedean local field,
that is a non-discrete totally disconnected locally compact
topological field. It is well known that $K$ is isomorphic either
to a finite extension of the field $\mathbb Q_p$ of $p$-adic
numbers (if $K$ has characteristic 0), or to the field of formal
Laurent series with coefficients from a finite field, if $K$ has
a positive characteristic. For a summary of main notions and results
regarding local fields see, for example, \cite{K2001}.

Any local field $K$ is endowed with an absolute value $|\cdot |_K$,
such that $|x|_K=0$ if and only if $x=0$, $|xy|_K=|x|_K\cdot |y|_K$,
$|x+y|_K\le \max (|x|_K,|y|_K)$. Denote $O=\{ x\in K:\ |x|_K\le 1\}$,
$P=\{ x\in K:\ |x|_K<1\}$. $O$ is a subring of $K$, and $P$ is an ideal
in $O$ containing such an element $\beta$ that
$P=\beta O$. The quotient ring $O/P$ is
actually a finite field; denote by $q$ its cardinality. We will
always assume that the absolute value is
normalized, that is $|\beta |_K=q^{-1}$. The normalized absolute
value takes the values $q^N$, $N\in \mathbb Z$. Note that for $K=\mathbb Q_p$
we have $\beta =p$ and $q=p$; the $p$-adic absolute value is normalized.

The additive group of any local field is self-dual, that is if
$\chi$ is a fixed non-constant complex-valued additive character of
$K$, then any other additive character can be written as
$\chi_a(x)=\chi (ax)$, $x\in K$, for some $a\in K$. Below we assume that $\chi$ is a rank
zero character, that is $\chi (x)\equiv 1$ for $x\in O$, while
there exists such an element $x_0\in K$ that $|x_0|_K=q$ and $\chi
(x_0)\ne 1$.

The above duality is used in the definition of the Fourier
transform over $K$. Denoting by $dx$ the Haar measure on the
additive group of $K$ (normalized in such a way that the measure
of $O$ equals 1) we write
$$
\widetilde{f}(\xi )=\int\limits_K\chi (x\xi )f(x)\,dx,\quad \xi
\in K,
$$
where $f$ is a complex-valued function from $L_1(K)$. As usual, the Fourier
transform $\mathcal F$ can be extended from $L_1(K)\cap L_2(K)$ to a
unitary operator on $L_2(K)$. If $\mathcal F f=\widetilde{f}\in L_1(K)$, we
have the inversion formula
$$
f(x)=\int\limits_K\chi (-x\xi )\widetilde{f}(\xi )\,d\xi .
$$

Working with functions on $K$ and operators upon them we often use standard integration formulas; see \cite{K2001,VVZ}. The simplest of them are as follows:

\begin{equation*}
\int\limits_{|x|_K\le q^n}dx=q^n;\quad \int\limits_{|x|_K=q^n}dx=\left( 1-\frac1q \right)q^n.
\end{equation*}

\begin{equation*}
\int\limits_{|x|_K\le q^n}|x|_K^{\alpha -1}\,dx=\frac{1-q^{-1}}{1-q^{-\alpha }}q^{\alpha n};\quad \text{here and above $n\in \mathbb Z,\alpha >0$}.
\end{equation*}

A function $f:\ K\to \mathbb C$ is said to be locally constant, if there exists such an integer $l$ that for any $x\in K$
$$
f(x+x')=f(x), \quad \text{whenever $|x'|\le q^{-l}$}.
$$
The smallest number $l$ with this property is called the exponent
of local constancy of the function $f$.

Let $\mathcal D(K)$ be the set of all locally constant
functions with compact supports; it is a vector space
over $\mathbb C$ with the topology of double inductive limit
$$
\mathcal D(K)=\varinjlim_{N\to \infty }\varinjlim_{l\to \infty }\mathcal D_N^l
$$
where $\mathcal D_N^l$ is the finite-dimensional space of functions supported in the ball
$B_N=\big\{ x\in K:$ $|x|\le q^N\big\}$ and having the exponents of local constancy $\le l$.
The strong conjugate space $\mathcal D'(K)$ is called the space of
Bruhat-Schwartz distributions.

The Fourier transform preserves the space $\mathcal D(K)$. Therefore the Fourier transform of a distribution defined
by duality acts continuously on $\mathcal D'(K)$. As in the case of $\mathbb R^n$, there exists a
well-developed theory of distributions over local fields; it includes
such topics as convolution, direct product, homogeneous distributions etc (see \cite{AKS,K2001,VVZ}).

\medskip
{\bf 2.2. Vladimirov's operator}. On a test function $\varphi \in \mathcal D(K)$, the fractional differentiation operator $D^\alpha$, $\alpha >0$, is defined as
\begin{equation}
\label{2.1}
\left( D^\alpha \varphi \right) (x)=\mathcal F^{-1}\left[ |\xi |_K^\alpha
(\mathcal F (\varphi ))(\xi )\right] (x).
\end{equation}
Note that $\D$ does not preserve $\mathcal D(K)$; see \cite{AKS} regarding the spaces of test functions and distributions preserved by this operator.

The operator $D^\alpha$ can also be represented as a hypersingular integral operator:
\begin{equation}
\label{2.2}
\left( D^\alpha \varphi \right) (x)=\frac{1-q^\alpha }{1-q^{-\alpha
-1}}\int\limits_K |y|_K^{-\alpha -1}[\varphi (x-y)-\varphi (x)]\,dy
\end{equation}
\cite{K2001,VVZ}.
In contrast to (\ref{2.1}), the expression in the right of (\ref{2.2}) makes sense for wider classes of functions. In particular, $\D$ is defined on constant functions and annihilates them. Denote for brevity $\theta_\alpha =\dfrac{1-q^\alpha }{1-q^{-\alpha -1}}$.

Below we consider the operator $\D$ on a radial function $u=u(|x|_K)$; here we identify the function $x\mapsto u(|x|_K)$ on $K$ with the function $|x|_K\mapsto u(|x|_K)$ on $q^{\mathbb Z}$. This abuse of notation does not lead to confusion.

The explicit expression of $\D u$ for a radial function $u$ satisfying some growth restrictions near the origin and infinity was found in \cite{K2014}. If $u=u(|x|_K)$ is such that
\begin{equation}
\label{2.3}
\sum\limits_{k=-\infty}^m q^k\left| u(q^k)\right| <\infty ,\quad \sum\limits_{l=m}^\infty q^{-\alpha l}\left| u(q^l)\right| <\infty,
\end{equation}
for some $m\in \mathbb Z$, then for each $n\in \mathbb Z$ the expression in the right-hand side of (\ref{2.2}) with $\varphi (x)=u(|x|_K)$ exists for $|x|_K=q^n$, depends only on $|x|_K$, and
\begin{multline}
\label{2.4}
(D^\alpha u)(q^n)=\theta_\alpha \left(1-\frac1q \right)q^{-(\alpha +1)n}\sum\limits_{k=-\infty}^{n-1} q^ku(q^k) +q^{-\alpha n-1}\frac{q^\alpha +q-2}{1-q^{-\alpha -1}}u(q^n)\\
+\theta_\alpha \left(1-\frac1q \right)\sum\limits_{l=n+1}^\infty q^{-\alpha l}u(q^l).
\end{multline}

Under the conditions (\ref{2.3}), the expression (\ref{2.4}) agrees also with the definition of $\D$ in terms of Bruhat-Schwartz distributions (see Chapter 2 of \cite{VVZ}).

\medskip
{\bf 2.3. The regularized integral}. The fractional integral mentioned in Introduction, was defined in \cite{K2014} initially for $\varphi \in \mathcal D(K)$ as follows:
\begin{equation}
\label{2.5}
(I^\alpha \varphi )(x)=\frac{1-q^{-\alpha}}{1-q^{\alpha -1}}\int\limits_{|y|_K\le |x|_K}\left( |x-y|_K^{\alpha -1}-|y|_K^{\alpha -1}\right) \varphi (y)\,dy,\quad \alpha \ne 1,
\end{equation}
and
\begin{equation}
\label{2.6}
(I^1\varphi )(x)=\frac{1-q}{q\log q}\int\limits_{|y|_K\le |x|_K}\left( \log |x-y|_K-\log |y|_K\right) \varphi (y)\,dy.
\end{equation}
Note that the integrals are taken, for each fixed $x\in K$, over bounded sets, and $(I^\alpha \varphi )(0)=0$. These properties are different from those of the anti-derivatives $D^{-\alpha }$ studied in \cite{VVZ}.

Let $u=u(|x|_K)$ be a radial function, such that
$$
\sum\limits_{k=-\infty}^m \max \left( q^k,q^{\alpha k}\right) \left| u(q^k)\right| <\infty ,\quad \text{if $\alpha \ne 1$},
$$
and
$$
\sum\limits_{k=-\infty}^m |k|q^k \left| u(q^k)\right| <\infty ,\quad \text{if $\alpha =1$},
$$
for some $m\in \mathbb Z$. Then \cite{K2014} $I^\alpha u$ exists, it is a radial function, and for any $x\ne 0$,
\begin{equation}
\label{2.7}
(I^\alpha u)(|x|_K)=q^{-\alpha}|x|_K^\alpha u(|x|_K)+\frac{1-q^{-\alpha}}{1-q^{\alpha -1}}\int\limits_{|y|_K< |x|_K}\left( |x|_K^{\alpha -1}-|y|_K^{\alpha -1}\right) u(|y|_K)\,dy,\quad \alpha \ne 1,
\end{equation}
and
\begin{equation}
\label{2.8}
(I^1 u)(|x|_K)=q^{-1}|x|_K u(|x|_K)+\frac{1-q}{q\log q}\int\limits_{|y|_K<|x|_K}\left( \log |x|_K-\log |y|_K\right) u(|y|_K)\,dy.
\end{equation}

\medskip
Another result from \cite{K2014} shows that $\II$ is indeed a right inverse of $\D$ on an appropriate class of radial functions. Namely,
suppose that for some $m\in \mathbb Z$,
\begin{equation}
\label{2.9}
\sum\limits_{k=-\infty}^m \max \left( q^k,q^{\alpha k}\right) \left| v(q^k)\right| <\infty ,\quad \sum\limits_{l=m}^\infty \left| v(q^l)\right| <\infty ,
\end{equation}
if $\alpha \ne 1$, and
\begin{equation}
\label{2.10}
\sum\limits_{k=-\infty}^m |k|q^k \left| v(q^k)\right| <\infty ,\quad
\sum\limits_{l=m}^\infty l\left| v(q^l)\right| <\infty ,
\end{equation}
if $\alpha =1$. Then there exists $\left( D^\alpha I^\alpha v\right) (|x|_K)=v(|x|_K)$ for any $x\ne 0$.

Note that the decay conditions at infinity in (\ref{2.9}) and (\ref{2.10}) cannot be dropped. This follows from the important identity
\begin{equation}
\label{2.11}
\II 1=0
\end{equation}
proved in \cite{K2014}.

Let us prove that under some decay conditions near the origin and infinity $\II$ is also a left inverse of $\D$. First we need an estimate for $\D u$. Below $C$ will denote various positive constants.

\medskip
\begin{prop}
Suppose that $u(0)=0$,
\begin{equation}
\label{2.12}
|u(q^n)|\le Cq^{dn},\quad n\le 0;
\end{equation}
\begin{equation}
\label{2.13}
|u(q^n)|\le Cq^{hn},\quad n\ge 0,
\end{equation}
where $d>\max (0,\alpha -1)$, $0\le h<\alpha$, and $h<\alpha -1$, if $\alpha >1$. Then the function $w=\D u$ satisfies, for any $m\in \mathbb Z$, the inequalities
\begin{equation}
\label{2.14}
\sum\limits_{k=-\infty}^m \max (q^k,q^{\alpha k})\left| w(q^k)\right| <\infty,\quad \sum\limits_{l=m}^\infty \left| w(q^l)\right| <\infty,
\end{equation}
if $\alpha \ne 1$, or
\begin{equation}
\label{2.15}
\sum\limits_{k=-\infty}^m |k|\cdot q^k\left| w(q^k)\right| <\infty,\quad \sum\limits_{l=m}^\infty |l|\cdot \left| w(q^l)\right| <\infty,
\end{equation}
if $\alpha =1$. Moreover,
\begin{equation}
\label{2.16}
\II \D u=u.
\end{equation}
\end{prop}

\medskip
{\it Proof}. Let $\alpha \ne 1$. Under (\ref{2.12}) and (\ref{2.13}), the conditions (\ref{2.3}) are satisfied, and we may use the expression (\ref{2.4}). Denote the three terms in (\ref{2.4}) by $S_1^{(n)}$, $S_2^{(n)}$ and $S_3^{(n)}$ respectively. We have for $n\ge 1$ that
$$
\sum\limits_{k=-\infty}^{n-1} q^k\left| u(q^k)\right|\le C\sum\limits_{k=-\infty}^0q^{(d+1)k}+C\sum\limits_{k=1}^{n-1}q^{(h+1)k}\le C_1q^{(h+1)n},
$$
so that $\left| S_1^{(n)}\right| \le Cq^{(h-\alpha )n}$ and $\sum\limits_{n=m}^\infty \left| S_1^{(n)}\right| < \infty$.

If $n\le 0$, then $\left| S_1^{(n)}\right| \le Cq^{(d-\alpha )n}$, whence
$$
\sum\limits_{k=-\infty}^m q^{\alpha k}\left| S_1^{(k)}\right| \le C\sum\limits_{k=-\infty}^m q^{dk}<\infty .
$$

Similar estimates for $S_2^{(k)}$ are obvious.

For $n\ge 1$,
$$
\left| S_3^{(n)}\right| \le \sum\limits_{l=n+1}^\infty q^{(h-\alpha )l}=C_2q^{(h-\alpha )n},
$$
so that $\sum\limits_{n=m}^\infty  \left| S_3^{(n)}\right|<\infty$. Next, let $0<\alpha <1$, so that $\max (q^k,q^{\alpha k})=q^{\alpha k}$ for $k\le 0$. We have
\begin{multline*}
\sum\limits_{n=-\infty}^{-2}q^{\alpha n}\sum\limits_{l=n+1}^\infty q^{-\alpha l}\left| u(q^l)\right| =\sum\limits_{n=-\infty}^{-2}q^{\alpha n}\left( \sum\limits_{l=n+1}^{-1}+\sum\limits_{l=0}^\infty \right) q^{-\alpha l}\left| u(q^l)\right| \\
\le C\sum\limits_{n=-\infty}^{-2}q^{\alpha n}\left( \sum\limits_{l=n+1}^{-1}q^{(d-\alpha )l}+\sum\limits_{l=0}^\infty q^{(h-\alpha )l}\right) =C\sum\limits_{l=-\infty}^{-1}q^{(d-\alpha )l}\sum\limits_{n=-\infty}^{l-1}q^{\alpha n}+C_3\sum\limits_{n=-\infty}^{-2}q^{\alpha n}\\
=C_4\sum\limits_{l=-\infty}^{-1}q^{dl}+C_3\sum\limits_{n=-\infty}^{-2}q^{\alpha n}<\infty.
\end{multline*}
Therefore
$$
\sum\limits_{k=-\infty}^m q^{\alpha k}\left| S_3^{(k)}\right| \le \infty,
$$
which implies the first inequality in (\ref{2.14}).

If $\alpha >1$, then $\max (q^k,q^{\alpha k})=q^k$ for $k\le 0$, and we perform similar calculations assuming that $d>\alpha -1$.

Finally, the proofs of the inequalities (\ref{2.15}) are similar to the above ones, since the factor $|k|$ does not influence the convergence of series with estimates exponential in $k$.

Let us prove (\ref{2.16}). Denote $v=\II \D u$. This is a legitimate object, since $\D u$ satisfies the conditions, under which $\II$ can be applied. Moreover, $\D v=\D \II (\D u)=\D u$, with a possible exception of the origin. Therefore $\D (v-u)$ is concentrated at the origin, which means (see e.g. Theorem 1.9 from \cite{K2001} or Section 6.3 in \cite{VVZ}) that $\D (v-u)=a\delta$, $a\in \mathbb R$. By Theorem 1, Section 9.3 of \cite{VVZ},
\begin{equation}
\label{2.17}
v-u=aD^{-\alpha}\delta +c
\end{equation}
where $c\in \RR$ and (\cite{AKS}, Example 9.2.1)
$$
\left( D^{-\alpha}\delta \right) (x)=\frac{|x|_K^{\alpha -1}}{\Gamma_K(\alpha )},\quad \Gamma_K(\alpha )=\frac{1-q^{\alpha -1}}{1-q^{-\alpha}},\text{ if $\alpha \ne 1$},
$$
$$
\left( D^{-1}\delta \right) (x)=\frac{q-1}{q\log q}\log |x|_K, \text{ if $\alpha =1$}.
$$
We need to prove that $a=c=0$.

Let $0<\alpha <1$. We saw above that $\left| S_1^{(n)}\right|, \left| S_2^{(n)}\right| \le Cq^{(d-\alpha )n}$, $n\le 0$. We also have that
$$
\left| S_3^{(n)}\right| \le C\sum\limits_{l=n+1}^\infty q^{-\alpha l}\left| u(q^l)\right| \le C_5+C_6\sum\limits_{l=n+1}^0q^{(d-\alpha )l}\le C_7q^{\min(0,d-\alpha )\cdot n},
$$
so that
\begin{equation}
\label{2.18}
\left| w(q^n)\right| \le Cq^{\min (0,d-\alpha )\cdot n},\quad n\le 0.
\end{equation}

Substituting (\ref{2.18}) into the expression (\ref{2.7}) for $\II w$ we find by a straightforward calculation that
\begin{equation}
\label{2.19}
\left| (\II w) (|x|_K)\right| \le C|x|_K^{\min (d,\alpha )} \to 0,\quad \text{as $x\to 0$}.
\end{equation}
Now it suffices to compare the assumption (\ref{2.12}), the estimate (\ref{2.19}) for $v=\II w$, and the relation (\ref{2.17}), to come to the conclusion that $a=c=0$, so that $v=u$.

The same argument works for $\alpha =1$. In this case the inequality (\ref{2.18}) remains valid, so that again $(I^1w)(q^n)\to 0$, as $n\to -\infty$, while $\log q^n\to -\infty$, as $n\to -\infty$.

Let $\alpha >1$. It follows from the above estimates of $S_1^{(n)},S_2^{(n)},S_3^{(n)}$ for $n\ge 1$ that
$$
\left| \left( \D u\right) (|x|_K)\right| \le C|x|_K^{h-\alpha},\quad |x|_K\ge 1.
$$
This implies the inequality
$$
|v(|x|_K)|\le C|x|_K^h, \quad |x|_K\ge 1
$$
(see the remark after Lemma 2 in \cite{K2014}). Comparing this with (\ref{2.17}) and our assumption (\ref{2.13}) we find that $a=0$. Then we consider as above the behavior of $u,v$ near the origin and find that $c=0$. $\qquad \blacksquare$

\medskip
\begin{cor}
Let $v=v_0+u$ where $v_0$ is a constant, $u$ satisfies the conditions of the above Proposition. Then $\II \D v=v-v_0$.
\end{cor}

\medskip
{\it Proof}. We have $\D v=\D u$, $\II \D v=u=v-v_0$. $\qquad \blacksquare$

\medskip
\section{The Cauchy problem}

{\bf 3.1. Local solvability.} Let us consider the problem (\ref{1.1})-(\ref{1.2}), that is
\begin{equation}
\label{3.1}
\left( \D u\right) (|t|_K) =f(|t|_K,u(|t|_K)),\quad 0\ne t\in K,
\end{equation}
with the initial condition
\begin{equation}
\label{3.2}
u(0)=u_0
\end{equation}
where the function $f:\ q^{\mathbb Z}\times \RR \to \RR$ satisfies the conditions
\begin{equation}
\label{3.3}
|f(|t|_K,x)|\le M;
\end{equation}
\begin{equation}
\label{3.4}
|f(|t|_K,x)-f(|t|_K,y)|\le F|x-y|,
\end{equation}
for all $t\in K,x,y\in \RR$.

With the problem (\ref{3.1})-(\ref{3.2}) we associate the integral equation
\begin{equation}
\label{3.5}
u(|t|_K)=u_0+\II f(|\cdot |_K,u(|\cdot |_K))(|t|_K).
\end{equation}
Note that, by the definition of $\II$, in order to compute $(\II \varphi )(|t|_K)$ for $|t|_K\le q^m$ ($m\in \mathbb Z$), one needs to know the function $\varphi$ in the same ball  $|t|_K\le q^m$. Therefore local solutions of the equation (\ref{3.5}) make sense, in contrast to solutions of (\ref{3.1}).

We will call a solution $u$ of (\ref{3.5}), if it exists, {\it a mild solution} of the Cauchy problem (\ref{3.1})-(\ref{3.2}). By the above Corollary, a solution $u$ of (\ref{3.1})-(\ref{3.2}), such that $u-u_0$ satisfies the conditions of Proposition, is a mild solution.

\medskip
\begin{teo}
Under the assumptions (\ref{3.3}),(\ref{3.4}), the problem (\ref{3.1})-(\ref{3.2}) has a unique local mild solution, that is the integral equation (\ref{3.5}) has a solution $u(|t|_K)$ defined for $|t|_K\le q^N$ where $N\in \mathbb Z$ is sufficiently small, and another solution $\overline{u}(|t|_K)$, if it exists, coincides with $u$ for $|t|_K\le q^{\overline{N}}$ where $\overline{N}\le N$.
\end{teo}

\medskip
{\it Proof}. We look for a solution of the equation (\ref{3.5}) as a limit of the sequence $\{ u_k\}$ where $u_0$ is the initial value from (\ref{3.5}),
\begin{equation}
\label{3.6}
u_k(|t|_K)=u_0+\II f(|\cdot |_K,u_{k-1}(|\cdot |_K))(|t|_K).
\end{equation}

We will use the following inequality proved in \cite{K2014}. Let
$$
I_{\alpha ,m}=\int\limits_{|y|_K<|t|_K}\left| |t|_K^{\alpha -1}-|y|_K^{\alpha -1}\right| |y|_K^{\alpha m}\,dy,\quad \alpha \ne 1,
$$
and
$$
I_{1,m}=\int\limits_{|y|_K<|t|_K}\left( \log |t|_K-\log |y|_K\right)|y|_K^m\,dy.
$$
Then
\begin{equation}
\label{3.7}
I_{\alpha ,m}=d_{\alpha ,m}|t|_K^{\alpha (m+1)},\quad 0<d_{\alpha,m}\le Aq^{-\alpha m},
\end{equation}
where $A>0$ does not depend on $m$.

Note in particular the case $m=0$, which implies the following inequality: if $\varphi$ is a bounded continuous function, then
\begin{equation}
\label{3.8}
\left| (\II \varphi )(|t|_K)\right| \le C\mu |t|_K^\alpha ,\quad \mu =\sup\limits_{t\in K}|\varphi (|t|_K)|.
\end{equation}
The meaning of this property is different from its classical counterparts, in view of the identity (\ref{2.11}).

Similarly, if $|\varphi (|s|_K)|\le \mu |s|_K^{\alpha m}$, then
\begin{equation}
\label{3.9}
\left| (\II \varphi )(|t|_K)\right| \le C\mu |t|_K^{\alpha (m+1)},\quad m=1,2,\ldots .
\end{equation}

Using (\ref{3.3}),(\ref{3.4}),(\ref{3.7}) and (\ref{3.8}),(\ref{3.9}) we find that
$$
|u_1(|t|_K)-u_0|\le CM|t|_K^\alpha ,
$$
$$
|u_2(|t|_K)-u_1(|t|_K)|=\left| \II [f(|\cdot|_K,u_1(|\cdot |_K))-f(|\cdot|_K,u_0)](|t|_K)\right|
$$
where
$$
|f(|s|_K,u_1(|s|_K))-f(|s|_K,u_0)|\le F|u_1(|s|_K)-u_0|\le CMF|s|_K^\alpha ,
$$
so that
$$
|u_2(|t|_K)-u_1(|t|_K)|\le C^2MF|t|_K^{2\alpha}.
$$

By induction, we find that
$$
|u_{k+1}(|t|_K)-u_k(|t|_K)|\le C^{k+1}MF^k|t|_K^{(k+1)\alpha}.
$$
This shows that the sequence $\{ u_k\}$ converges uniformly on the ball $\{ |t|_K\le T\}$, if $T$ is sufficiently small, to a limit $u(|t|_K)$. By (\ref{2.5}) and (\ref{2.6}), we obtain that
$$
\II f(|\cdot|_K,u_k(|\cdot |_K))(|t|_K)\longrightarrow f(|\cdot|_K,u(|\cdot |_K))(|t|_K),\quad |t|_K\le T,
$$
so that $u(|t|_K)$ is indeed a local solution of the equation (\ref{3.5}).

Suppose we have another local solution $\uu (|t|_K)$. We have by (\ref{3.8}) that
$$
|u(|t|_K)-\uu (|t|_K)|\le 2MC|t|_K^\alpha.
$$
On the other hand,
$$
u(|t|_K)-\uu (|t|_K)=\II f(|\cdot|_K,u(|\cdot |_K))(|t|_K)-\II f(|\cdot|_K,\uu (|\cdot |_K))(|t|_K),
$$
and by (\ref{3.4}) and (\ref{3.9}),
$$
|u(|t|_K)-\uu (|t|_K)|\le 2MC^2F|t|_K^{2\alpha}.
$$

Repeating these arguments we obtain by induction that
$$
|u(|t|_K)-\uu (|t|_K)|\le 2MC^m F^{m-1}|t|_K^{m\alpha}.
$$
For sufficiently small $|t|_K$, this iteration with $m\to \infty$ implies the equality $u(|t|_K)=\uu (|t|_K)$. $\qquad \blacksquare$

\bigskip
{\bf 3.2. Extension of solutions.} Let us study the possibility to continue the local solution constructed in Theorem 1 to a solution of the integral equation (\ref{3.5}) defined for all $t\in K$.

Suppose that the conditions of Theorem 1 are satisfied, and we obtained a local solution $u(|t|_K)$, $|t|_K\le q^N$, $N\in \mathbb Z$. Let $\alpha \ne 1$. In order to find a solution for $|t|_K=q^{N+1}$, we have, by virtue of (\ref{2.7}), to solve the equation
\begin{equation}
\label{3.10}
u(q^{N+1})=u_0+v_0^{(N)}+q^{\alpha N}f(q^{N+!},u(q^{N+1}))
\end{equation}
where
\begin{equation}
\label{3.11}
v_0^{(N)}=\frac{1-q^{-\alpha}}{1-q^{\alpha -1}}\int\limits_{|y|_K\le q^N}\left( q^{(N+1)(\alpha -1)}-|y|_K^{\alpha -1}\right) f(|y|_K,u(|y|_K))\,dy.
\end{equation}
is a known constant. A similar equation can be written for $\alpha =1$. Then the above procedure, if it is successful, is repeated for all $l>N$.

The following result is an immediate consequence of Banach's fixed point theorem.

\medskip
\begin{teo}
Suppose that the conditions of Theorem 1 are satisfied, as well as the following Lipschitz condition:
\begin{equation}
\label{3.12}
|f(q^l,x)-f(q^l,y)|\le F_l |x-y|,\quad x,y\in \RR ,\ l\in \mathbb Z,
\end{equation}
where $0<F_l<q^{-\alpha l}$ for each $l\in \mathbb Z$. Then a local solution of the equation (\ref{3.5}) admits a continuation to a global solution defined for all $t\in K$.
\end{teo}

\bigskip
{\bf 3.3. From an integral equation to a differential one.} Let us study conditions, under which the above continuation procedure leads to a solution of the problem (\ref{3.1})-(\ref{3.2}). As before, we assume the conditions (\ref{3.3}),(\ref{3.4}) and (\ref{3.12}). In addition, we will assume that
\begin{equation}
\label{3.13}
|f(q^l,x)|\le Cq^{-\beta l},\quad l\ge 1,\quad \text{for all $x\in \RR$,}
\end{equation}
where $\beta >\alpha$.

Let us estimate the mild solution $u(|t|_K)$ obtained by iterations for $|t|_K\le q^N$, and then extended like in (\ref{3.10}). This solution is automatically continuous at the origin. For each $l\ge N$, the solution satisfies the equation similar to (\ref{3.10}):
$$
u(q^{l+1})=u_0+v_0^{(l)}+q^{\alpha l}f(q^{l+1},u(q^{l+1})).
$$

For $l\ge 1$, let us write $v_0^{(l)}=v_{0,1}^{(l)}+v_{0,2}^{(l)}$ where the summands correspond to the integration over the sets $\{ |y|_K\le 1\}$ and $\{ |y|_K\ge q\}$ respectively. Then
$$
\left| v_{0,1}^{(l)}\right| \le C\left[ q^{(l+1)(\alpha -1)}+1\right] ;
$$
\begin{multline*}
\left| v_{0,2}^{(l)}\right| \le Cq^{(l+1)(\alpha -1)}\int\limits_{q\le |y|_K\le q^l}|y|_K^{-\beta}dy+C\int\limits_{q\le |y|_K\le q^l}|y|_K^{\alpha -1-\beta}dy\\
=C_1q^{(l+1)(\alpha -1)}\sum\limits_{j=1}^l q^{-\beta j+j}+C_1\sum\limits_{j=1}^l q^{j(\alpha -1-\beta)+j}\le C_2\left( 1+q^{(\alpha -\beta)l}\right) .
\end{multline*}

It follows from (\ref{3.10}) and (\ref{3.13}) that the conditions (\ref{2.3}) of the existence of $\D u$ are satisfied. In addition, it follows from (\ref{3.13}) that $f(|t|_K,u(|t|_K))$ satisfies the condition (\ref{2.9}) (or (\ref{2.10}), for $\alpha =1$), under which $\D \II f=f$. We come to the following result.

\medskip
\begin{teo}
Under the assumptions (\ref{3.3}),(\ref{3.4}), (\ref{3.12}) and (\ref{3.13}), the mild solution obtained by the iteration process with subsequent continuation, satisfies the equation (\ref{3.1}).
\end{teo}

\section*{Acknowledgments}
This work was funded in part under the budget program of Ukraine No. 6541230 ``Support to the development of priority research trends'' and under the research work "Markov evolutions in real and p-adic spaces" of the Dragomanov National Pedagogical  University of Ukraine.

\end{document}